\documentclass[11pt,leqno]{article}
\usepackage{amsmath, amsthm,  amsfonts}
\input{amssym.def}

\numberwithin{equation}{section}

\def\and{\hbox{\quad and \quad}}

\catcode`\@=12
\newtheorem{theorem}{Theorem}

\theoremstyle{remark}

\newtheorem{remark}{Remark}

\input{epsf.tex}

\begin{document}

\vskip 0.3cm
\begin{center}
{\LARGE  {\bf Corresponding 
constant mean curvature surfaces in hyperbolic and 
Euclidean 3-spaces}}
\end{center}

\begin{center}
{\bf Wayne Rossman, Magdalena Toda}
\end{center}

\begin{center}
{\bf Dedicated to the memory of Hongyou Wu}
\end{center}

{\abstract We make observations about constant mean curvature surfaces in Euclidean $3$-space 
and their dual surfaces, and the resulting pairs of 
surfaces in hyperbolic $3$-space under the Lawson 
correspondence.} \vskip 0.4cm



\vskip 0.4cm
It is well-known that any constant mean curvature 
(CMC) surface in Euclidean $3$-space $R^3$ 
can be parameterized in conformal curvature 
line coordinates away from umbilics, and such a parameterization 
is called {\it isothermic}.  
A characteristic property of an isothermic immersion is the 
(local) existence of a {\it dual} surface $f^d$ 
(via Christoffel transform, see \cite{Chris} and 
\cite{HerPe}).  Like in \cite{HerPe}, 
we say that two surfaces in Euclidean $3$-space form a 
Christoffel pair if they induce 
conformally equivalent metrics and have parallel tangent planes 
with opposite orientations. 
The Christoffel transform is an involution, and each 
of the surfaces of a Christoffel pair is called a Christoffel 
transform, or dual, of the other.  The dual 
(Christoffel transform) is unique up to homothety and translation. 

Let $\Sigma$ be a simply-connected Riemann surface and let 
$f: \Sigma \to R^3$ be an isothermic
constant mean curvature (CMC) $H$ immersion that is not the round sphere.
Assume $H \neq 0$, so the surface is not minimal, 
and let $N: \Sigma \to S^2$ ($S^2$ is the round sphere of radius $1$ 
centered at the origin in $R^3$) be the unit normal vector to $f$.  
An appropriately scaled and positioned Christoffel 
transform of $f$ will be its parallel constant mean curvature surface 
$f^d=f+H^{-1} N: \Sigma \to R^3$ with normal 
$N^d = - N$ (see, for example, \cite{HerPe}), which has the same 
constant mean curvature as $f$ itself, i.e. $H_d=H$.  
It is well known that this Christoffel transform is also a Darboux 
transform; in fact, it was shown in \cite{HerPe} 
that the Christoffel transform $f^d$ of an isothermic surface $f$ is also a 
Darboux transform of $f$ if and only if $f$ has non-zero constant 
mean curvature $H$. 

Taking the isothermic CMC immersions $f$ and $f^d$ as above, we can 
rescale the complex coordinate $z$ of $\Sigma$ so that $f$ and $f^d$ have 
induced conformal metrics that are inverse to each other, and we then call 
$z$ a normalized isothermic coordinate.  We find that 
\[ H = H_d = 2 Q = 2 Q_d \; , \] 
where $Q$ and $Q_d$ are the Hopf differential functions of $f$ and $f^d$, 
respectively.  
The metric, Hopf differential function and mean curvature of $f$ are 
\[ I = e^{2u} (dx^2+dy^2) \; , \;\;\; Q = \langle f_{zz}, N \rangle \; , 
\;\;\; H = 2 e^{-2u} \langle f_{z\bar z} , N \rangle \; , \] 
where $\langle \cdot , \cdot \rangle$ is the bilinear extension to 
complex $3$-space $C^3$ of the usual Euclidean metric for $R^3$.  
Then $f^d$ has the corresponding data 
\[ I_d = e^{-2u} (dx^2+dy^2) \; , \;\;\; Q_d = \langle f^d_{zz}, N^d \rangle 
= Q \; , \;\;\; H_d = 2 e^{2u} \langle f^d_{z\bar z} , N^d \rangle = H \; . \] 

Let the hyperbolic $3$-space $H^3$ be given as an isometric 
submanifold of Minkowski $4$-space $R^{3,1}$ (with $(+++-)$ 
metric) via $H^3 = \{ (x_1,x_2,x_3,x_0) \in R^{3,1} \, | \, 
x_0 > 0, x_0^2-x_1^2-x_2^2-x_3^2=1 \}$, 
which can also be written in the Hermitean matrix model as 
\[ H^3 = \left\{ \left. \begin{pmatrix} x_0+x_3 & x_1-i x_2 \\ 
x_1+i x_2 & x_0-x_3 \end{pmatrix} \, \right| \, 
x_0 > 0, x_0^2-x_1^2-x_2^2-x_3^2=1 \right\}= \{ F \bar{F}^t \, | \, 
F \in \text{SL}_2C \} \; . \]  
With respect to the Hermitean matrix model, points in 
$R^{3,1}$ can be given as 
\[ X = \begin{pmatrix} x_0+x_3 & x_1-i x_2 \\ 
x_1+i x_2 & x_0-x_3 \end{pmatrix} \] (where $\det X$ is 
not necessarily $1$), and then the metric for $R^{3,1}$ is 
\[ \langle X,Y \rangle = - \tfrac{1}{2} \text{tr}\left(
X \begin{pmatrix} 0 & -i \\ i & 0 \end{pmatrix} Y^t 
\begin{pmatrix} 0 & -i \\ i & 0 \end{pmatrix} \right) \]
for $X, Y \in R^{3,1}$.  

We will be considering homotheties $s f$, resp. $s f^d$, 
of the immersion $f$, resp. $f^d$, for 
$s \in R \setminus \{ 0 \}$.  Then, the Lawson correspondent 
$f_1$, resp. $f_1^d$, in $H^3$, of 
$s f$, resp. $s f^d$, has the geometric data 
\begin{equation}\label{eqn:star1} 
I_1 = s^2 e^{2u} (dx^2+dy^2) \; , \;\;\; Q_1 = s Q \; , 
\;\;\; H_1 = \sqrt{(s^{-1} H)^2+1} \; , \;\;\; 
\text{resp.} \end{equation} 
\begin{equation}\label{eqn:star1pt5} 
I_{d,1} = s^2 e^{-2u} (dx^2+dy^2) \; , \;\;\; Q_{d,1} 
= s Q \; , \;\;\; H_{d,1} = \sqrt{(s^{-1} H)^2+1} \; . \end{equation}

Before stating our result, we note that the 
extended $r$-unitary frame for $f$ can be given by 
a solution $F$ to 
the following Lax system (see \cite{Ai2}, \cite{Bo}, \cite{DoPeWu}, 
\cite{KKRS2}, \cite{To}) 
\[ F_z = F U \; , \;\;\; F_{\bar z} = F V \; , 
\]\[ 
U = \frac{1}{2} \begin{pmatrix}
-u_z & 2 e^{-u} \lambda^{-1} Q \\ - H e^u & u_z 
\end{pmatrix} \; , \;\;\; 
V = \frac{1}{2} \begin{pmatrix}
u_{\bar z} & H e^u \\ - 2 e^{-u} \lambda Q & -u_{\bar z} 
\end{pmatrix} \; , \] 
whose compatibility condition is the Gauss equation 
\[ 4 u_{z \bar z} - 4 Q^2 e^{-2 u} + H^2 e^{2u} = 0 \; . \] 

For the Euclidean immersion $f$, the parameter $\lambda$ 
belongs to the unit circle $S^1$, and then the frame $F$ 
belongs to the $\text{SU}_2$-valued loop group.  The actual frame 
of $f$ is represented by $F|_{\lambda = 1}$.  
However, to create an $r$-unitary frame, we fix $r$ as some 
positive real number less than one, and $\lambda$ 
is taken to be in the open annular domain in the complex plane 
between the circles of radius $r$ and radius $1/r$ centered at 
the origin. For more details on $r$-unitary frames, 
see \cite{KKRS2}.  Let $D$ be the diagonal matrix 
\[ D = \begin{pmatrix} 1/\sqrt{\lambda} & 0 \\ 0 & \sqrt{\lambda}
\end{pmatrix} \; . \] 
Now we can state our result:

\begin{theorem}
Let $f$ be a CMC immersion in $R^3$ as above, without 
umbilic points and with normalized isothermic coordinate $z$, and 
with extended $r$-unitary frame $F$.  
Choose a value of $\lambda$ so that $r < \lambda < 1$.  
Then the two surfaces $F \bar{F}^t$ and 
$F D \overline{F D}^t$ in $H^3$, each evaluated at that value of 
$\lambda$, are both isothermic CMC immersions 
in $H^3$, and the following hold: 
\begin{enumerate}
\item $F \bar{F}^t$ is the Lawson correspondent to a homothety 
       of $f^d$; 
\item $F D \overline{F D}^t$ is the Lawson 
       correspondent to a homothety of $f$; 
\item $F D \overline{F D}^t$ is an equidistant (parallel) surface 
       to $F \bar{F}^t$, with mean curvature 
       opposite in sign to that of $F \bar{F}^t$. 
\end{enumerate}
\end{theorem}

\begin{proof}
The metric, Hopf differential function and (hyperbolic) mean curvature 
of $F \bar{F}^t$, resp. $F D \overline{F D}^t$, are 
\begin{equation}\label{eqn:star2} 
Q^2 e^{-2u} (\lambda-\lambda^{-1})^2 
(dx^2+dy^2) \; , \;\;\; 
\tfrac{1}{2} Q H (\lambda^{-1}-\lambda) \; , \;\;\; 
\frac{\lambda^{-1}+\lambda}{\lambda^{-1}-\lambda} \; , \;\;\; 
\text{resp.} 
\end{equation}  
\begin{equation}\label{eqn:star3} 
Q^2 e^{2u} (\lambda-\lambda^{-1})^2 (dx^2+dy^2) \; , \;\;\; 
\tfrac{1}{2} Q H (\lambda-\lambda^{-1}) \; , \;\;\; 
\frac{\lambda+\lambda^{-1}}{\lambda-\lambda^{-1}} \; . 
\end{equation}  
These surfaces are isothermic CMC immersions in $H^3$.  

The metric, Hopf differential function and mean curvature 
of the Lawson correspondent of the homothety $s f^d$ of $f^d$, 
resp. $s f$ of $f$, are as in \eqref{eqn:star1pt5}, resp. as in 
\eqref{eqn:star1}, and that 
data will be equal to the data \eqref{eqn:star2} 
for $F \bar{F}^t$, resp. the data \eqref{eqn:star3} 
for $F D \overline{F D}^t$, if 
$s = \tfrac{1}{2} H (\lambda^{-1}-\lambda)$, 
resp. $s = \tfrac{1}{2} H (\lambda-\lambda^{-1})$, proving the 
first and second itemized statements of the theorem, as well as 
the second half of the third item.  

Write $\lambda = e^q$ for some value $q < 0$.  
Noting that the normal vector to $(F \bar{F}^t)|_{\lambda = e^q}$ is 
\[ \hat N = \left. \left( F \begin{pmatrix} 
1 & 0 \\ 0 & -1 \end{pmatrix} \bar{F}^t \right) 
\right|_{\lambda = e^q} \; , \] we have that 
$(F D \overline{F D}^t)|_{\lambda=e^q} = 
(\cosh q) (F \bar{F}^t)|_{\lambda=e^q} - (\sinh q) \hat N$.
Thus the distance between the surfaces 
$F \bar{F}^t$ and $F D \overline{F D}^t$ is $-q$, completing 
the proof of the third item.  
\end{proof}

\begin{remark}
The surface $F \bar{F}^t$ was used in \cite{KKRS2} 
to construct CMC trinoids in $H^3$, thus those trinoids in $H^3$ are not 
Lawson correspondents to the CMC trinoids in $R^3$ with frame $F$, 
but rather to the duals of those CMC trinoids in $R^3$.  This shift 
to the dual surface 
in order to close periods on non-simply-connected CMC surfaces in $H^3$ 
occured in \cite{RoUmYa} as well.  
\end{remark}

\begin{remark}
One simple example to which the theorem can be applied is round cylinders 
in $H^3$, given by a round cylinder in $R^3$ with data 
\[ u=0 \; , \;\;\; H = \tfrac{1}{2} \; , \;\;\; 
Q = \tfrac{1}{4} \]
and 
\[ F = \begin{pmatrix} 
\cosh \gamma & (\sqrt{\lambda})^{-1} \sinh \gamma \\ 
\sqrt{\lambda} \sinh \gamma & \cosh \gamma 
\end{pmatrix} \; , \;\;\; 
\gamma = \frac{1}{4} i \left( \frac{z}{\sqrt{\lambda}} + 
\sqrt{\lambda} \bar z \right) \; . 
\]
\end{remark}

{\bf Acknowledgement.}  The authors thank Udo Hertrich-Jeromin 
for noticing an error in their initial computations.  They also note that 
Udo Hertrich-Jeromin found a different proof of the result here, using 
linear conserved quantities instead of unitary frames.

\bigskip

Magdalena Toda,

Department of Mathematics and Statistics,

Texas Tech University,

Lubbock, Texas 79409-1042,

U.S.A.

{\em magda.toda@ttu.edu}

\bigskip

\

\

Wayne Rossman,

Department of Mathematics, Faculty of Science, 

Kobe University,

Rokko, Kobe 657-8501,

Japan

{\em wayne@math.kobe-u.ac.jp}

\end{document}